\documentclass[english]{article}
\usepackage[latin1]{inputenc}
\usepackage{amsmath,amsthm,amssymb,babel}
\newtheorem{teo}{Theorem}

\newtheorem{prop}{Proposition}

\newtheorem{lemma}{Lemma}

\def\proof{{\it Proof.}\ }

\def\eq#1{(\ref{#1})}

\def\neweq#1{\begin{equation}\label{#1}}
\def\endeq{\end{equation}}

\def\phi{\varphi}
\def\RR{{\mathbb R} }

\def\di{\displaystyle}
\def\ri{\rightarrow}

\date{}

\title{\sc Multiple solutions for a nonlinear and non-homogeneous problem in Orlicz-Sobolev spaces}
\author{\sc   Mihai Mih\u ailescu$\,^{a}$ and Du\v san Repov\v s$\,^{b}$ \\
 \small $^a\,$Department of Mathematics, Central European
University,  1051 Budapest, Hungary\\
\small
$^b\,$Faculty of Mathematics and Physics, and Faculty of Education, University of Ljubljana,\\
\small
 POB 2964, Ljubljana, Slovenia 1001\\
\small
E-mail addresses:  {\tt mmihailes@yahoo.com}\qquad {\tt dusan.repovs@guest.arnes.si} }

\begin{document}
\maketitle
\noindent{\small{\sc Abstract}. We study a non-homogeneous  boundary value problem in a smooth bounded domain in
$\RR^N$. We prove the existence of at least two nonnegative and non-trivial weak solutions. Our approach
relies on Orlicz-Sobolev spaces theory combined with adequate variational methods and a variant of Mountain Pass Lemma.  
}

\section{Introduction and preliminary results}
Let $\Omega$ be a bounded domain in $\RR^N$ ($N\geq 3$) with smooth boundary $\partial\Omega$. Assume that  $a:(0,\infty)\ri\RR$ is a function
such that the mapping $\phi:\RR\rightarrow\RR$, defined by
$$\phi(t)=\left\{\begin{array}{lll}
a(|t|)t, &\mbox{for}&
t\neq 0\\
0, &\mbox{for}& t=0\,,
\end{array}\right.$$
is an odd, increasing homeomorphisms from $\RR$ onto $\RR$. This paper studies a nonlinear
boundary value problem of the type
\begin{equation}\label{1}
\left\{\begin{array}{lll}
-{\rm div}(a(|\nabla u|)\nabla u)=\lambda f(x,u), &\mbox{for}&
x\in\Omega\\
u=0, &\mbox{for}& x\in\partial\Omega
\end{array}\right.
\end{equation}
where $f:\Omega\times\RR\rightarrow\RR$ is a Carath\'eodory function and $\lambda$ is a positive parameter.

In order to go further we introduce the functional space setting where problem \eq{1} will be discussed. In this context
we note that the operator in the divergence form is not homogeneous and thus, we introduce an Orlicz-Sobolev space setting
for problems of this type.

The first general existence results using the theory of monotone
operators in Orlicz-Sobolev spaces were obtained by Donaldson
\cite{dona1} and Gossez \cite{G}. Other recent works that put the
problem into this framework are contained in Cl\'ement {\it et
al.} \cite{Clem1,Clem2}, Garc\'ia-Huidobro {\it et al.}
\cite{Gar}, Gossez and Man\`asevich \cite{gosman}, Le and Schmitt
\cite{les}, etc. In these papers, the existence results are
obtained by means of variational techniques, monotone operator
methods, or fixed point and degree theory arguments. The goal of
our paper is to present a new multiplicity result for equations
involving nonhomogeneous operators. Thus, it supplements the
aforementioned results in the aspect that most of the papers
guarantee existence but not multiplicity of solutions.

We start by recalling some basic facts about Orlicz spaces. For more
details we refer to the books by D. R. Adams and L. L. Hedberg
\cite{AHed}, R. Adams \cite{A} and M. M. Rao and Z. D. Ren
\cite{rao} and the papers by Ph. Cl\'ement et al. \cite{Clem1,
Clem2}, M. Garc\'ia-Huidobro et al. \cite{Gar} and J. P. Gossez
\cite{G}.

For  $\phi:\RR\rightarrow\RR$ introduced at the start of the
paper, we define
$$\Phi(t)=\int_0^t\phi(s)\;ds,\;\;\;\Phi^\star(t)=\int_0^t\phi^{-1}(s)\;ds,\qquad \mbox{for all}\ t\in\RR\,.$$
We observe that $\Phi$ is a {\it Young function}, that is,
$\Phi(0)=0$, $\Phi$ is convex, and
$\lim_{x\ri\infty}\Phi(x)=+\infty$. Furthermore, since $\Phi(x)=0$
if and only if $x=0$, $\lim_{x\ri 0}\Phi(x)/x=0$, and $\lim_{x\ri
\infty}\Phi(x)/x=+\infty$, then $\Phi$ is called an {\it
$N$--function}. The function $\Phi^\star$ is called the {\it
complementary} function of $\Phi$, and it satisfies
$$\Phi^\star (t)=\sup\{st-\Phi(s);\ s\geq 0\},\qquad\mbox{for all $t\geq 0$}\,.$$
 We also observe that $\Phi^\star$ is also an  $N$--function and the Young's inequality holds
 $$st\leq\Phi(s)+\Phi^\star (t),\qquad\mbox{for all $s,t\geq 0$}\,.$$

The Orlicz space $L_{\Phi}(\Omega)$ defined by the $N$--function $\Phi$
(see \cite{AHed,A,Clem1}) is the space of measurable functions $u:\Omega\ri\RR$ such that
$$\|u\|_{L_{\Phi}}:=\sup\left\{\int_\Omega uv\;dx;\ \int_\Omega(\Phi)^\star (|g|)\;dx\leq 1\right\}<\infty\,.$$
Then $(L_{\Phi}(\Omega),\|\,\cdot\, \|_{L_{\Phi}} )$ is a Banach space whose norm is equivalent to the
 Luxemburg norm
$$\|u\|_{\Phi} :=\inf\left\{k>0;\ \int_\Omega\Phi\left(\frac{u(x)}{k}
\right)\;dx\leq 1\right\}.$$ For Orlicz spaces  H\"older's
inequality reads as follows (see \cite[Inequality 4, p. 79]{rao}):
$$\int_\Omega uvdx\leq 2\,\|u\|_{L_{\Phi}}\, \|v\|_{L_{(\Phi)^\star}}\qquad\mbox{for all $u\in L_{\Phi}(\Omega)$
and $v\in L_{(\Phi)^\star}(\Omega)$}\,.$$

Next, we introduce the Orlicz-Sobolev spaces. We denote by $W^1L_{\Phi}(\Omega)$ the Orlicz-Sobolev space defined
by
$$W^1L_{\Phi}(\Omega):=\left\{u\in L_{\Phi}(\Omega);\;\frac{\partial u}
{\partial x_i}\in L_{\Phi}(\Omega),\;i=1,...,N\right\}\,.$$
This is a Banach space with respect to the norm
$$\|u\|_{1,\Phi}:=\|u\|_{\Phi}+\||\nabla u|\|_{\Phi}\,.$$
We also define the Orlicz-Sobolev space $W_0^1L_{\Phi}(\Omega)$ as
the closure of $C_0^\infty(\Omega)$ in $W^1L_{\Phi}(\Omega)$. By
\cite[Lemma 5.7]{G} we may consider on $W_0^1L_{\Phi}(\Omega)$
the equivalent norm
$$\|u\|:=\||\nabla u|\|_{\Phi}.$$
For an easier manipulation of the spaces defined above, we define
$$\phi_0:=\inf_{t>0}\frac{t\phi(t)}{\Phi(t)}\;\;{\rm and}\;\;\phi^0:=\sup_{t>0}\frac{t\phi(t)}{\Phi(t)}\,.$$
In this paper we assume that we have
\begin{equation}\label{acdc0}1<\phi_0\leq\frac{t\phi(t)}{\Phi(t)}\leq\phi^0<\infty,\;\;\;\forall\;t\geq
0\,.\end{equation} The above relation implies that  $\Phi$ satisfies the $\Delta_2$-condition, i.e.
\begin{equation}\label{acdc}
\Phi(2t)\leq K\Phi(t),\;\;\;\forall\;t\geq 0\,,
\end{equation}
where $K$ is a positive constant (see  \cite[Proposition
2.3]{fourier}).

On the other hand, the following relations hold
\begin{equation}\label{zz}
\small
\|u\|^{\phi^0}\leq\int_\Omega\Phi(|\nabla u|)\;dx\leq\|u\|^{\phi_0},\;\;\;\forall\;u\in W_0^1L_{\Phi}(\Omega),\;\|u\|<1\,,
\end{equation}
\begin{equation}\label{zz1}
\small
\|u\|^{\phi_0}\leq\int_\Omega\Phi(|\nabla u|)\;dx\leq\|u\|^{\phi^0},\;\;\;\forall\;u\in W_0^1L_{\Phi}(\Omega),\;\|u\|>1\,,
\end{equation}
(see, e.g. \cite[Lemma 1]{AA}).

Furthermore, in this paper we shall assume that  the function
$\Phi$ satisfies the following condition
\begin{equation}\label{acdc1}
{\rm the}\;{\rm function}\;[0,\infty)\ni t\rightarrow\Phi(\sqrt{t})\;{\rm is}\;{\rm convex}\,.
\end{equation}
Conditions \eq{acdc} and \eq{acdc1} assure that  the Orlicz space
$L_{\Phi}(\Omega)$ is a uniformly convex space and thus, a
reflexive Banach space (see \cite[Proposition 2.2]{fourier}). This
fact implies that also the Orlicz-Sobolev space
$W_0^1L_{\Phi}(\Omega)$ is a reflexive Banach space.
\smallskip

\noindent{\bf Remark.} We point out certain examples of functions $\phi:\RR\rightarrow\RR$ which are odd, increasing homeomorphisms from $\RR$
onto $\RR$ and satisfy conditions \eq{acdc0} and \eq{acdc1}. For more details the reader can consult \cite[Examples 1-3, p. 243]{Clem2}.

1) Let
$$\phi(t)=p|t|^{p-2}t,\;\;\;\forall\;t\in\RR\,,$$
with $p>1$. For this function it can be proved that
$$\phi_0=\phi^0=p\,.$$
Furthermore, in this particular case the corresponding Orlicz space $L_\Phi(\Omega)$ is the classical Lebesgue space
$L^p(\Omega)$ while the Orlicz-Sobolev space $W_0^1L_\Phi(\Omega)$ is the classical Sobolev space $W_0^{1,p}(\Omega)$.
We will use the classical notations to denote the Orlicz-Sobolev spaces in this particular case.

2) Consider
$$\phi(t)=\log(1+|t|^s)|t|^{p-2}t,\;\;\;\forall\;t\in\RR\,,$$
with $p$, $s>1$. In this case it can be proved that
$$\phi_0=p,\;\;\;\;\;\phi^0=p+s\,.$$

3) Let
$$\phi(t)=\frac{|t|^{p-2}t}{\log(1+|t|)},\;\;\;{\rm if}\;t\neq 0,\;\;\;\phi(0)=0\,,$$
with $p>2$. In this case we have
$$\phi_0=p-1,\;\;\;\;\;\phi^0=p\,.$$
\smallskip

 \section{The main result}
In this paper we study problem \eq{1} in the special case when
$$f(x,t)=t^{p-1}-t^{q-1}\,,$$
with
\begin{equation}\label{steluta}
1<q<p<\phi_0
\end{equation}
and $t\geq 0$.

More precisely, we consider the degenerate boundary value problem
\begin{equation}\label{2}
\left\{\begin{array}{lll}
-{\rm div}(a(|\nabla u|)\nabla u)=\lambda(u^{p-1}-u^{q-1}), &\mbox{for}& x\in\Omega\\
u=0, &\mbox{for}& x\in\partial\Omega\\
u\geq 0,&\mbox{for}& x\in\Omega.
\end{array}\right.
\end{equation}

We say that $u\in W_0^1L_{\Phi}(\Omega)$ is a {\it weak solution} of problem \eq{2} if $u\geq 0$ a. e. in $\Omega$ and
$$\int_\Omega a(|\nabla u|)\nabla u\cdot\nabla v\;dx-\lambda\int_\Omega u^{p-1}v\;dx+\lambda\int_\Omega
u^{q-1}v\;dx=0$$
for all $v\in W_0^1L_{\Phi}(\Omega)$.
\smallskip

Our main result asserts that problem \eq{2} has at least two nontrivial weak solutions provided that $\lambda>0$ is large enough. More precisely,
we prove

\begin{teo}\label{t1}
Assume that condition \eq{steluta} is fulfilled and
\begin{equation}\label{163276}\phi^0<\min\left\{N,\frac{N\phi_0}{N-\phi_0}\right\}\,.\end{equation}
Then there exists $\lambda^\star>0$ such that for all $\lambda>\lambda^\star$ problem \eq{2} has at least two distinct non-negative, nontrivial
weak solutions.
\end{teo}

\noindent{\bf Remark.} We point out that our result was inspired
by \cite[Theorem 1.2]{P}, where a related property was proved in
the case of the $p$-Laplace operators. The extension from
$p$-Laplace operator to the differential operators involved in
\eq{2}  is not trivial, since the new operators have a more
complicated structure than the $p$-Laplace operator, for example
they are non-homogeneous.
\smallskip

Finally, we mention that a similar study regarding the existence and multiplicity of solutions for equations involving the $p(x)$-Laplace operator
can be found in Mih\u ailescu and R\u adulescu \cite{RoyalSoc}.

\section{Proof of Theorem \ref{t1}}
Let $E$ denote the generalized Sobolev
space $W_0^1L_{\Phi}(\Omega)$.

Define the energy functional
$I:E\rightarrow\RR$ by
$$I(u)=\int_\Omega \Phi(|\nabla u|)\;dx-\frac{\lambda}
{\gamma}\int_\Omega u_+^p\;dx+\frac{\lambda}{\beta}\int_\Omega
u_+^q\;dx\,,$$
where $u_+(x)=\max\{u(x),0\}$.

We remember that $u\in E$ implies $u_+$, $u_-\in E$ and
$$\nabla u_+=\left\{\begin{array}{lll}
0, &\mbox{if}& [u\leq 0]\\
\nabla u, &\mbox{if}& [u>0],
\end{array}\right. \qquad
\nabla u_-=\left\{\begin{array}{lll}
0, &\mbox{if}& [u\geq 0]\\
\nabla u, &\mbox{if}& [u<0]
\end{array}\right.$$
where $u_\pm=\max\{\pm u(x),0\}$ for all $x\in\Omega$ (see, e.g.
page 52 in \cite{Clem1}). That fact and some standards arguments
assure that  functional $I$ is well-defined on $E$ and $I\in
C^1(E,\RR)$ with the derivative given by
$$\langle I^{'}(u),v\rangle=\int_\Omega a(|\nabla u|)\nabla u\cdot\nabla v\;dx-\lambda\int_\Omega u_+^{p-1}v\;dx
+\lambda\int_\Omega u_+^{q-1}v\;dx,$$
for all $u$, $v\in E$.

\noindent{\bf Remark.} We point out that if $u$ is a critical point of $I$  then using the above information and condition \eq{acdc0} we have
\begin{eqnarray*}
0&=&\langle I^{'}(u),u_-\rangle=\int_\Omega a(|\nabla u|)\nabla u
\cdot\nabla u_-\;dx-\lambda\int_\Omega(u_+)^{p-1}u_-\;dx
+\lambda\int_\Omega(u_+)^{q-1}u_-\;dx\\
&=&\int_\Omega a|\nabla u|)\nabla u\cdot\nabla u_-\;dx
=\int_\Omega a(|\nabla u_-|)|\nabla u_-|^2\;dx
\geq\phi_0\int_\Omega\Phi(|\nabla u_-|)\;dx.
\end{eqnarray*}
By the above estimates and relation \eq{zz} we deduce that $u\geq 0$. It follows that the nontrivial critical points of $I$ are non-negative
solutions of \eq{2}.
\smallskip

The above remark shows that we can prove Theorem \ref{t1} using the critical point theory. More exactly, we first show that for $\lambda>0$ large
enough, the functional $I$ has a global minimizer $u_1\geq 0$ such that $I(u_1)<0$. Next, by means of the Mountain Pass Theorem, a second critical
point $u_2$ with $I(u_2)>0$ is obtained.

\begin{lemma}\label{le4}
There exists $\lambda_1>0$ such that
$$\lambda_1=\inf\limits_{u\in E,\;\|u\|>1}\frac{\di\int_\Omega
\Phi(|\nabla u|)\;dx}{\di\int_\Omega|u|^{\phi_0}\;dx}\,.$$
\end{lemma}
\proof First, we note that by condition \eq{acdc0} we can deduce
that $E$ is continuously embedded in the classical Sobolev space
$W_0^{1,\phi_0}(\Omega)$. Consequently, $E$ is continuously
embedded in the classical Lebesgue space $L^{\phi_0}(\Omega)$. It
follows that there exists $C>0$ such that
$$\|u\|\geq C\|u\|_{L^{\phi_0}(\Omega)},\;\;\;\forall\; u\in E.$$
On the other hand, by \eq{zz1} we have
$$\int_\Omega\Phi(|\nabla u|)\;dx\geq\|u\|^{\phi_0},\;\;\;\forall\;
u\in E\;{\rm with}\;\|u\|>1.$$
Combining the above inequalities we obtain
$$\int_\Omega\Phi(|\nabla u|)\;dx\geq C^{\phi_0}\int_\Omega|u|^{\phi_0}\;dx,\;\;\;\forall\; u\in E\;{\rm with}\;
\|u\|>1.$$
The proof of Lemma \ref{le4} is complete.\qed

\begin{prop}\label{p2}
(i) The functional $I$ is bounded from below and coercive.\\
(ii) The functional $I$ is weakly lower semi-continuous.
\end{prop}
\proof
(i) Since $1<q<p<\phi_0$ we have
$$\lim\limits_{t\rightarrow\infty}\frac{\di\frac{1}{p}t^p-
\di\frac{1}{q}t^q}{t^{\phi_0}}=0.$$
Then for any $\lambda>0$ there exists $C_\lambda>0$ such that
$$\lambda\left(\frac{1}{p}t^p-\frac{1}{q}t^q
\right)\leq\frac{\lambda_1}{2}t^{\phi_0}+C_\lambda,\;\;\;\forall\;
t\geq 0,$$ where $\lambda_1$ was defined in Lemma \ref{le4}.

The above inequality and condition \eq{zz1} show that for any $u\in E$ with $\|u\|>1$ we have
\begin{eqnarray*}
I(u)&\geq&\int_\Omega\Phi(|\nabla u|)\;dx-
\frac{\lambda_1}{2}\int_\Omega|u|^{\phi_0}\;dx-C_\lambda\mu(\Omega)\\
&\geq&\frac{1}{2}\int_\Omega\Phi(|\nabla u|)\;dx-
C_\lambda\mu(\Omega)\\
&\geq&\frac{1}{2}\|u\|^{\phi_0}-C_\lambda\mu(\Omega).
\end{eqnarray*}
This shows that $I$ is bounded from below and coercive.
\smallskip

\noindent (ii) Similar arguments as those used in the proof of \cite[Theorem 2]{jmaa} (see also \cite[Lemma 4.3]{fourier}) show that the
functional $I_0:E\rightarrow\RR$ defined by
\begin{equation}\label{081577}I_0(u)=\int_\Omega\Phi(|\nabla u|)\;dx\,,\end{equation}
is weakly lower semi-continuous. We justify  that $I$ is weakly lower semi-continuous. Let
$(u_n)\subset E$ be a sequence which converges weakly to $u$ in $E$.
Since $I_0$ is weakly lower semi-continuous we have
\begin{equation}\label{3}
I_0(u)\leq\liminf\limits_{n\rightarrow\infty}I_0(u_n).
\end{equation}
On the other hand, since $E$ is compactly embedded in $L^p(\Omega)$ and $L^q(\Omega)$ it follows that $({u_n}_+)$ converges strongly to $u_+$ both
in $L^p(\Omega)$ and in $L^q(\Omega)$. (The compact embedding of $E$ into $L^p(\Omega)$ and $L^q(\Omega)$ is a direct consequence
of the fact that $E$ is continuously embedded in the classical Sobolev space $W_0^{1,\phi_0}(\Omega)$ combined with condition \eq{steluta}.)
This fact together with relation \eq{3} imply
$$I(u)\leq\liminf\limits_{n\rightarrow\infty}I(u_n)\,.$$
Therefore, $I$ is weakly lower semi-continuous.
The proof of Proposition \ref{p2} is complete.  \qed

From Proposition \ref{p2} and Theorem 1.2 in \cite{S} we deduce
that there exists $u_1\in E$ a global minimizer of $I$. The
following result implies that $u_1\neq 0$, provided that $\lambda$
is sufficiently large.

\begin{prop}\label{p3}
There exists $\lambda^\star>0$ such that $\inf_E I<0$.
\end{prop}

\proof
Let $\Omega_1\subset\Omega$ be a compact subset, large enough and
$u_0\in E$ be such that $u_0(x)=t_0$ in $\Omega_1$ and $0\leq u_0(x)
\leq t_0$ in $\Omega\setminus\Omega_1$, where $t_0>1$ is chosen
such that
$$\frac{1}{p}t_0^p-\frac{1}{q}t_0^q>0.$$
We have
\begin{eqnarray*}
\frac{1}{p}\int_\Omega u_0^p\;dx-\frac{1}{q}\int_\Omega
u_0^q\;dx&\geq&\frac{1}{p}\int_{\Omega_1}u_0^p\;dx-
\frac{1}{q}\int_{\Omega_1}u_0^q\;dx-\frac{1}{q}
\int_{\Omega\setminus\Omega_1}u_0^q\;dx\\
&\geq&\frac{1}{p}
\int_{\Omega_1}u_0^p\;dx-\frac{1}{q}\int_{\Omega_1}u_0^q
\;dx-\frac{1}{q}\;t_0^q\;\mu(\Omega\setminus\Omega_1)>0
\end{eqnarray*}
and thus $I(u_0)<0$ for $\lambda>0$ large enough.
The proof of Proposition \ref{p3} is complete.  \qed

Since Proposition \ref{p3} holds  it follows that $u_1\in E$ is a
nontrivial weak solution of problem \eq{2}.
\bigskip

Fix $\lambda\geq\lambda^\star$. Set
$$
g(x,t)=\left\{\begin{array}{lll}
0, &\mbox{for}& t<0\\
t^{p-1}-t^{q-1}, &\mbox{for}& 0\leq t\leq u_1(x)\\
u_1(x)^{p-1}-u_1(x)^{q-1}, &\mbox{for}& t>u_1(x)
\end{array}\right.
$$
and
$$G(x,t)=\int_0^tg(x,s)\;ds.$$
Define the functional $J:E\rightarrow\RR$ by
$$J(u)=\int_\Omega \Phi(|\nabla u|)\;dx-\lambda \int_\Omega G(x,u)
\;dx.$$
The same arguments as those used for functional $I$ imply that
$J\in C^1(E,\RR)$ and
$$\langle J^{'}(u),v\rangle=\int_\Omega a(|\nabla u|)\nabla u
\cdot\nabla v\;dx-\lambda\int_\Omega g(x,u)v\;dx,$$
for all $u$, $v\in E$.

On the other hand, we point out that if $u\in E$ is a
critical point of $J$ then $u\geq 0$. The proof can be carried out
as in the case of functional $I$.

Next, we prove
\begin{lemma}\label{l5}
If $u$ is a critical point of $J$ then $u\leq u_1$.
\end{lemma}
\proof
We have
\begin{eqnarray*}
0&=&\langle J^{'}(u)-I^{'}(u_1),(u-u_1)_+\rangle\\
&=&\int_\Omega(a(|\nabla u|)\nabla u-a(|\nabla u_1|)\nabla u_1)
\cdot\nabla(u-u_1)_+\;dx-\lambda\int_\Omega[g(x,u)-(u_1^{p-1}-u_1^{q-1})](u-u_1)_+\;dx\\
&=&\int_{[u>u_1]}(a(|\nabla u|)\nabla u-a(|\nabla u_1|)\nabla u_1)
\cdot\nabla(u-u_1)\;dx\,.
\end{eqnarray*}
Notice that since $\phi$ is increasing in $\RR$ we have for each
$\xi$ and $\psi\in\RR^N$
$$(\phi(|\xi|)-\phi(|\psi|))(|\xi|-|\psi|)\geq 0\,,$$
with equality if and only if $\xi=\psi$. Thus, we can deduce that
$$(a(|\xi|)|\xi|-a(|\psi|)|\psi|)(|\xi|-|\psi|)\geq 0\,,$$
for all  $\xi,\psi\in\RR^N$, with equality if and only if $\xi=\psi$. On the other hand, some simple computations show that
$$(a(|\xi|)\xi-a(|\psi|)\psi)\cdot(\xi-\psi)\geq (a(|\xi|)|\xi|-\phi(|\psi|)|\psi|)(|\xi|-|\psi|)\,,$$
for all  $\xi,\psi\in\RR^N$. Consequently, we conclude that
$$(a(|\xi|)\xi-a(|\psi|)\psi)\cdot(\xi-\psi)\geq 0\,,$$
for all  $\xi,\psi\in\RR^N$, with equality if and only if $\xi=\psi$.

Using the above pieces of information we deduce that the above equality holds if and only if $\nabla u=\nabla u_1$. It follows that $\nabla u(x)=
\nabla u_1(x)$ for all $x\in\omega:=\{y\in\Omega;\;\;u(y)>u_1(y)\}$.
Hence
$$\int_\omega\Phi(|\nabla(u-u_1)|)\;dx=0$$
and thus
$$\int_\Omega\Phi(|\nabla(u-u_1)_+|)\;dx=0.$$
By relation \eq{zz} we obtain
$$\|(u-u_1)_+\|=0.$$
We obtain that $(u-u_1)_+=0$ in $\Omega$, that is, $u\leq u_1$ in
$\Omega$. The proof of Lemma \ref{l5} is complete.  \qed

In the following we determine a critical point $u_2\in E$ of $J$
such that $J(u_2)>0$ via the Mountain Pass Theorem. By the above lemma
we will deduce that $0\leq u_2\leq u_1$ in $\Omega$. Therefore
$$g(x,u_2)=u_2^{p-1}-u_2^{q-1}\;\;\;{\rm and}\;\;\;
G(x,u_2)=\frac{1}{p}u_2^p-\frac{1}{q}u_2^q$$
and thus
$$J(u_2)=I(u_2)\;\;\;{\rm and}\;\;\;J^{'}(u_2)=I^{'}(u_2).$$
More precisely,  we find
$$I(u_2)>0=I(0)>I(u_1)\qquad\mbox{and}\qquad I^{'}(u_2)=0\,.$$
This shows that $u_2$ is a weak solution of problem \eq{2} such
that $0\leq u_2\leq u_1$, $u_2\neq 0$ and $u_2\neq u_1$.

In order to find $u_2$ described above we prove
\begin{lemma}\label{l6}
There exists $\rho\in(0,\|u_1\|)$ and $a>0$ such that
$J(u)\geq a$, for all $u\in E$ with $\|u\|=\rho.$
\end{lemma}

\proof
Let $u\in E$ be fixed, such that $\|u\|<1$. It is clear that there
exists $\delta>1$ such that
$$\frac{1}{p}t^p-\frac{1}{q}t^q\leq 0,\;\;\;
\forall t\in[0,\delta].$$
For $\delta$ given above we define
$$\Omega_u=\{x\in\Omega;\;u(x)>\delta\}.$$
If $x\in\Omega\setminus\Omega_u$ with $u(x)<u_1(x)$ we have
$$G(x,u)=\frac{1}{p}u_+^p-\frac{1}{q}u_+^q\leq 0.$$
If $x\in\Omega\setminus\Omega_u$ with $u(x)>u_1(x)$ then
$u_1(x)\leq\delta$ and we have
$$G(x,u)=\frac{1}{p}u_1^p-\frac{1}{q}u_1^q\leq 0.$$
Thus we deduce that
$$G(x,u)\leq 0,\;\;\;{\rm on}\;\Omega\setminus\Omega_u.$$
Provided that $\|u\|<1$ by  relation \eq{zz} we get
\begin{equation}\label{s1}
\begin{array}{lll}
J(u)&\geq&\di\int_\Omega\di\Phi(|\nabla u|)\;dx-
\lambda\di\int_{\Omega_u}G(x,u)\;dx\\
&\geq&\|u\|^{\phi^0}-\lambda\di\int_{\Omega_u}G(x,u)\;dx
\end{array}
\end{equation}
By relation \eq{163276} it follows that
$\phi^0<\phi_0^\star:=\frac{N\phi_0}{N-\phi_0}$. On the other
hand, as we already pointed out, by condition \eq{acdc0} we deduce
that $E$ is continuously embedded in the classical Sobolev space
$W_0^{1,\phi_0}(\Omega)$. Consequently, there exists
$s\in(\phi^0,\frac{N\phi_0}{N-\phi_0})$ such that $E$ is
continuously embedded in the classical Lebesgue space
$L^{s}(\Omega)$. Thus, there exists a positive constant $C>0$ such
that
$$\|u\|_{L^s(\Omega)}\leq C\|u\|,\;\;\;\forall u\in E.$$
Using the definition of $G$, H\"older's inequality and the above estimate, we obtain
\begin{equation}\label{s2}
\begin{array}{lll}
\lambda\di\int_{\Omega_u} G(x,u)\;dx&=&\lambda\di\int_{\Omega_u\cap
[u<u_1]}\left(\di\frac{1}{p}u_+^p-\di\frac{1}{q}
u_+^q\right)\;dx+\lambda\di\int_{\Omega_u\cap[u>u_1]}\left(
\di\frac{1}{p}u_1^p-\di\frac{1}{q}u_1^q\right)
\;dx\\
&\leq&\di\frac{2\lambda}{p}\di\int_{\Omega_u}u_+^p\;dx\\
&\leq&\di\frac{2\lambda}{p}\di\int_{\Omega_u}u_+^{\phi^0}\;dx\\
&\leq&\di\frac{2\lambda}{p}\left(\di\int_{\Omega_u}u_+^{s}
\;dx\right)^{\phi^0/s}[\mu(\Omega_u)]^{1-\phi^0/q}\\
&\leq&C\di\frac{2\lambda}{p}[\mu(\Omega_u)]^{1-\phi^0/s}\|u\|^{\phi^0}.
\end{array}
\end{equation}
By \eq{s1} and \eq{s2} we infer that it is enough to show that
$\mu(\Omega_u)\rightarrow 0$ as $\|u\|\rightarrow 0$ in order to
prove Lemma \ref{l6}.

Let $\epsilon>0$. We choose $\Omega_\epsilon\subset\Omega$ a compact
subset, such that $\mu(\Omega\setminus{\Omega_\epsilon})<\epsilon$.
We denote by $\Omega_{u,\epsilon}:=\Omega_u\cap\Omega_\epsilon$. Then
it is clear that
$$C^{\phi^0}\|u\|^{\phi^0}\geq\left(\int_\Omega|u|^{s}
\;dx\right)^{\phi^0/s}\geq\left(\int_{\Omega_{u,\epsilon}}|u|^{s}\;dx
\right)^{\phi^0/s}\geq\delta^{\phi^0}[\mu(\Omega_{u,\epsilon})]^{\phi^0/s}.$$
The above inequality implies that $\mu(\Omega_{u,\epsilon})
\rightarrow 0$ as $\|u\|\rightarrow 0$.

Since $\Omega_u\subset\Omega_{u,\epsilon}\cup(\Omega\setminus
\Omega_\epsilon)$ we have
$$\mu(\Omega_u)\leq\mu(\Omega_{u,\epsilon})+\epsilon$$
and $\epsilon>0$ is arbitrary. We find that $\mu(\Omega_u)\rightarrow
0$ as $\|u\|\rightarrow 0$.
This concludes the proof of Lemma \ref{l6}. \qed

\begin{lemma}\label{l7}
The functional $J$ is coercive.
\end{lemma}
\proof
For each $u\in E$ with $\|u\|>1$ by  relation \eq{zz1} and H\"older's inequality  we have
\begin{eqnarray*}
J(u)&\geq&\int_\Omega\Phi(|\nabla u|)\;dx-\lambda\int_{[u>u_1]}G(x,u)\;dx-\lambda\int_{[u<u_1]}G(x,u)\;dx\\
&\geq&\|u\|^{\phi_0}-\frac{\lambda}{p}\int_{[u>u_1]}u_1^p\;dx+\frac{\lambda}{q}\int_{[u>u_1]}u_1^q\;dx-
\frac{\lambda}{p}\int_{[u<u_1]}u_+^p\;dx+\frac{\lambda}{q}\int_{[u<u_1]}u_+^q\;dx\\
&\geq&\|u\|^{\phi_0}-\frac{\lambda}{p}\int_{\Omega}u_1^p\;dx-\frac{\lambda}{p}\int_{\Omega}u_+^p\;dx\\
&\geq&\|u\|^{\phi_0}-\frac{\lambda}{p}[\mu(\Omega)]^{1-p/\phi_0}C_1\|u\|^p-C_2\\
&\geq&\|u\|^{\phi_0}-C_23\|u\|^p-C_2,
\end{eqnarray*}
where $C_1$, $C_2$ and $C_3$ are positive constants. Since $p<\phi_0$ the above inequality implies that $J(u)\rightarrow\infty$ as
$\|u\|\rightarrow\infty$, that is, $J$ is coercive. The proof of Lemma \ref{l7} is complete.  \qed

The following result yields a sufficient condition which ensures that a
weakly convergent sequence
in $E$ converges strongly, too.

\begin{lemma}\label{l8}
Assume that the sequence $(u_n)$ converges weakly to $u$ in $E$ and
$$\limsup\limits_{n\rightarrow\infty}\int_\Omega a(|\nabla u_n|)\nabla u_n\cdot(\nabla u_n-\nabla u)\;dx\leq 0\,.$$
Then $(u_n)$ converges strongly to $u$ in $E$.
\end{lemma}

\proof
Since $u_n$ converges weakly to $u$ in $E$ implies that it follows that $(\|u_n\|)$ is a bounded sequence of real numbers. That fact and relations
\eq{zz} and \eq{zz1} imply thet the sequence $(I_0(u_n))$ is bounded, where $I_0$ is defined by relation \eq{081577}. Then, up to to a subsequence,
we deduce that $I_0(u_n)\rightarrow c$. Furthermore, the weak lower semi-continuity of $I_0$ (pointed out above) implies
$$I_0(u)\leq\liminf\limits_{n\rightarrow\infty}I_0(u_n)=c\,.$$
On the other hand, since $I_0$ is convex (because $\Phi$ is convex), we have
$$I_0(u)\geq I_0(u_n)+\langle I_0^{'}(u_n),u-u_n\rangle\,.$$
Next, by the hypothesis  $\limsup\limits_{n\rightarrow\infty}
\int_\Omega a(|\nabla u_n|)\nabla u_n\cdot(\nabla u_n-\nabla u)\;dx\leq 0$,
we conclude that $I_0(u)=c$.

Taking into account that $(u_n+u)/2$ converges weakly to $u$ in $E$ and using again the weak lower semi-continuity of $I_0$ we find
\begin{equation}\label{ecc6}
c=I_0(u)\leq\liminf\limits_{n\rightarrow\infty}I_0\left(\frac{u_n+u}{2}\right)\,.
\end{equation}
We assume by contradiction that $u_n$ does not converge to $u$ in $E$. Then by \eq{zz} it follows that there exist $\epsilon>0$ and a
subsequence $(u_{n_m})$ of $(u_n)$ such that
\begin{equation}\label{ecc7}
I_0\left(\frac{u_{n_m}-u}{2}\right)\geq\epsilon,\;\;\;\forall\; m\,.
\end{equation}
On the other hand, relations \eq{acdc} and \eq{acdc1} enable us to apply \cite[Lemma 2.1]{lamperti} in order to obtain
\begin{equation}\label{ecc8}
\frac{1}{2}I_0(u)+\frac{1}{2}I_0(u_{n_m})-I_0\left(\frac{u+u_{n_m}}{2}\right)\geq I_0\left(\frac{u-u_{n_m}}{2}\right)\geq\epsilon,
\;\;\;\forall\; m\,.
\end{equation}
Letting $m\rightarrow\infty$ in the above inequality we obtain
$$c-\epsilon\geq\limsup\limits_{m\rightarrow\infty}I_0\left(\frac{u+u_{n_m}}{2}\right)\;dx\,,$$
and that is a contradiction with \eq{ecc6}. It follows that $u_n$ converges strongly to $u$ in $E$ and Lemma \ref{l8} is proved.  \qed

\medskip
{\sc Proof of Theorem \ref{t1} completed.}
Using Lemma \ref{l6} and the Mountain Pass Theorem (see \cite{AR} with the variant given by Theorem 1.15 in \cite{W}) we deduce that there exists
a sequence $(u_n)\subset E$ such that
\begin{equation}\label{MPT}
J(u_n)\rightarrow c>0\;\;\;{\rm and}\;\;\;J^{'}(u_n)\rightarrow 0
\end{equation}
where
$$c=\inf\limits_{\gamma\in\Gamma}\max\limits_{t\in[0,1]}
J(\gamma(t))$$
and
$$\Gamma=\{\gamma\in C([0,1],E);\;\gamma(0)=0,\;\gamma(1)=u_1\}.$$
By relation \eq{MPT} and Lemma \ref{l7} we obtain that $(u_n)$ is bounded and thus passing eventually to a subsequence, still denoted by $(u_n)$,
we may assume that there exists $u_2\in E$ such that $u_n$ converges weakly to $u_2$. Since $E$ is compactly embedded in $L^i(\Omega)$ for any
$i\in[1,\phi_0]$, it follows that $u_n$ converges strongly to $u_2$ in $L^i(\Omega)$ for all $i\in[1,\phi_0]$. Hence
$$\langle I_0^{'}(u_n)-I_0^{'}(u_2),u_n-u_2\rangle=
\langle J^{'}(u_n)-J^{'}(u_2),u_n-u_2\rangle+\lambda\int_\Omega
[g(x,u_n)-g(x,u_2)](u_n-u_2)\;dx=o(1)\,,$$
as $n\ri\infty$, where $I_0$ is defined by relation \eq{081577}. By Lemma \ref{l8}  we deduce that $u_n$ converges strongly to $u_2$ in $E$
and using relation \eq{MPT} we find
$$J(u_2)=c>0\;\;\;{\rm and}\;\;\;J^{'}(u_2)=0.$$
Therefore, $J(u_2)=c>0$ and $J^{'}(u_2)=0$. By Lemma \ref{l5} we deduce that $0\leq u_2\leq u_1$ in $\Omega$. Therefore
$$g(x,u_2)=u_2^{p-1}-u_2^{q-1}\;\;\;{\rm and}\;\;\;
G(x,u_2)=\frac{1}{p}u_2^p-\frac{1}{q}u_2^q$$
and thus
$$J(u_2)=I(u_2)\;\;\;{\rm and}\;\;\;J^{'}(u_2)=I^{'}(u_2).$$
We conclude that $u_2$ is a critical point of $I$ and thus a solution of \eq{2}. Furthermore, $I(u_2)=c>0$ and $I(u_2)>0>I(u_1)$. Thus $u_2$ is
not trivial and $u_2\neq u_1$. The proof of Theorem \ref{t1} is now complete.  \qed
\bigskip

\noindent {\bf Acknowledgements.} This research was supported by
Slovenian Research Agency grants P1-0292-0101 and J1-2057-0101.

\end{document}